\documentclass[12pt,reqno]{amsart}
\usepackage{amsmath,amssymb}

\newtheorem{theorem}{Theorem}[section]
\newtheorem{lemma}[theorem]{Lemma}
\newtheorem{proposition}[theorem]{Proposition}
\newtheorem{corollary}[theorem]{Corollary}

\theoremstyle{definition}

\newtheorem{remark}[theorem]{Remark}

\numberwithin{equation}{section}

\begin{document}

\title[Equivariant vector bundles on projective line]{Equivariant vector bundles over
the complex projective line}

\author[I. Biswas]{Indranil Biswas}

\address{Department of Mathematics, Shiv Nadar University, NH91, Tehsil
Dadri, Greater Noida, Uttar Pradesh 201314, India}

\email{indranil.biswas@snu.edu.in, indranil29@gmail.com}

\author[F.-X. Machu]{Francois-Xavier Machu}

\address{ESIEA, 74 bis Av. Maurice Thorez, 94200 Ivry-sur-Seine, France}

\email{fx.machu@gmail.com}



\date{}

\begin{abstract}
Let $G$ be a finite abelian group acting faithfully on ${\mathbb C}{\mathbb P}^1$ via holomorphic
automorphisms. In \cite{DF2} the $G$--equivariant algebraic
vector bundles on $G$--invariant affine open subsets of ${\mathbb C}{\mathbb P}^1$
were classified. We classify the $G$--equivariant algebraic vector bundles on ${\mathbb C}{\mathbb P}^1$.
\end{abstract}

\maketitle

\section{Introduction}

This work was inspired by \cite{DF2}. The set-up of \cite{DF2} is the following. A finite abelian group
$G$ is acting faithfully on ${\mathbb C}{\mathbb P}^1$ via holomorphic automorphisms. Let $X\, \subset\,
{\mathbb C}{\mathbb P}^1$ be an affine subspace preserved by the action of $G$, so ${\mathbb C}{\mathbb P}^1
\setminus X$ is a nonempty finite subset preserved by the action of $G$. The $G$--equivariant algebraic
vector bundles on $X$ were classified in \cite{DF2}. Our aim here is to classify the
$G$--equivariant algebraic vector bundles on ${\mathbb C}{\mathbb P}^1$.

The group $G$ can be of two types. Either it is a cyclic group or it is isomorphic to
$({\mathbb Z}/2{\mathbb Z}) \oplus ({\mathbb Z}/2{\mathbb Z})$.

When $G$ is a finite cyclic group, the following are proved (see Theorem \ref{th1} and Remark
\ref{rem1}):

\begin{enumerate}

\item Any $G$--equivariant holomorphic vector bundle on ${\mathbb C}{\mathbb P}^1$ splits into a
direct sum of $G$--equivariant holomorphic line bundles.

\item Every holomorphic line bundle on ${\mathbb C}{\mathbb P}^1$ admits a $G$--equivariant structure.

\item The $G$--equivariant structures on any given holomorphic line bundle $L$ are parametrized
by the group of complex characters of $G$.
\end{enumerate}

When $G$ is isomorphic to $({\mathbb Z}/2{\mathbb Z}) \oplus ({\mathbb Z}/2{\mathbb Z})$, the following
are proved (see Lemma \ref{lem2}, Theorem \ref{thm2} and Remark \ref{rem2}):

\begin{enumerate}
\item The holomorphic line bundle ${\mathcal O}_{{\mathbb C}{\mathbb P}^1}(2n)$
admits a $G$--equivariant structure for every $n\, \in\, {\mathbb Z}$.
The $G$--equivariant structures on ${\mathcal O}_{{\mathbb C}{\mathbb P}^1}(2n)$ are parametrized
by the group of complex characters of $G$.

\item The holomorphic vector bundle ${\mathcal O}_{{\mathbb C}{\mathbb P}^1}(2n-1)
\oplus {\mathcal O}_{{\mathbb C}{\mathbb P}^1}(2n-1)$ admits a $G$--equivariant structure
for every $n\, \in\, {\mathbb Z}$. The $G$--equivariant structures on ${\mathcal O}_{{\mathbb C}{\mathbb P}^1}
(2n-1)\oplus {\mathcal O}_{{\mathbb C}{\mathbb P}^1}(2n-1)$ are described in Remark \ref{rem3}.

\item For any $n\, \in\, {\mathbb Z}$, the holomorphic line bundle ${\mathcal O}_{{\mathbb C}
{\mathbb P}^1}(2n-1)$ does not admit any $G$--equivariant structure.

\item Any $G$--equivariant vector bundle $E$ admits a holomorphic decomposition into $G$--equivariant
vector bundles
$$
E \,\,=\,\, \left(\bigoplus_{i=1}^M {\mathcal O}_{{\mathbb C}{\mathbb P}^1}(2m_i)\right)
\oplus \left(\bigoplus_{j=1}^N {\mathcal O}_{{\mathbb C}{\mathbb P}^1}(2n_i+1)^{\oplus 2}\right),
$$
where $m_i$ and $n_i$ are integers and $M+2N\,=\, {\rm rank}(E)$.
\end{enumerate}

\section{Equivariant vector bundles on projective line}\label{se2}

The standard action of $\text{GL}(2,{\mathbb C})$ on ${\mathbb C}^2$ produces an action
of $\text{PGL}(2,{\mathbb C})$ on the complex projective line ${\mathbb C}{\mathbb P}^1\,=\,
{\mathbb C}\cup \{\infty\}$. This way $\text{PGL}(2,{\mathbb C})$ gets identified with
the holomorphic automorphism group ${\rm Aut}({\mathbb C}{\mathbb P}^1)$ of
${\mathbb C}{\mathbb P}^1$. Let $G$ be a finite abelian subgroup of $\text{PGL}(2,{\mathbb C})$.
We know the following:
\begin{enumerate}
\item either $G$ is a finite cyclic group, or

\item $G\,=\, ({\mathbb Z}/2{\mathbb Z}) \oplus ({\mathbb Z}/2{\mathbb Z})$.
\end{enumerate}
(See \cite[Section 1.1]{Do}, \cite{Sp}.) If $G\,=\, {\mathbb Z}/n{\mathbb Z}$, then it is
conjugate to the subgroup of $\text{PGL}(2,{\mathbb C})$ generated by
$\begin{pmatrix}\alpha & 0\\ 0 & 1\end{pmatrix}$, where $\alpha$ is a primitive
$n$-th root of unity. In this case the two point $0$ and $\infty$ are fixed by the action of $G$
on ${\mathbb C}{\mathbb P}^1$. If $G\,=\, ({\mathbb Z}/2{\mathbb Z}) \oplus ({\mathbb Z}/2{\mathbb Z})$, then
it is conjugate to the subgroup of $\text{PGL}(2,{\mathbb C})$ generated by
\begin{equation}\label{e1}
A_1\, :=\, \begin{pmatrix}-1 & 0\\ 0 & 1\end{pmatrix}\ \ \,\text{ and }\ \ \,
A_2\,:=\, \begin{pmatrix}0 & 1\\ 1 & 0\end{pmatrix}.
\end{equation}
In this case, no element of ${\mathbb C}{\mathbb P}^1$ is fixed by the action of $G$ on ${\mathbb C}{\mathbb P}^1$.
The action of $A_1$ fixes the points $0$ and $\infty$; the action of $A_2$ fixes
$1$ and $-1$ while the action of $A_1A_2$ fixes $\sqrt{-1}$ and $-\sqrt{-1}$.

Let $\varpi\, :\, E\, \longrightarrow\, {\mathbb C}{\mathbb P}^1$ be a holomorphic vector bundle.
A $G$--equivariant structure on $E$ is a holomorphic action of the group $G$
\begin{equation}\label{e2}
\phi\,\, :\,\, G\times E \,\, \longrightarrow\,\, E
\end{equation}
on the total space of $E$ such that for every $\gamma\, \in\, G$, the map
\begin{equation}\label{e3}
\phi_\gamma\,\, :\,\, E \,\, \longrightarrow\,\, E
\end{equation}
defined by $z\,\longmapsto\, \phi(\gamma,\, z)$ is an automorphism of the
vector bundle $E$ over that automorphism $\gamma$ on ${\mathbb C}{\mathbb P}^1$. Therefore, we have
$\varpi\circ\phi_\gamma\,=\, \gamma\circ \varpi$ for all $\gamma\, \in\, G$.
Also, for any $x\, \in\, {\mathbb C}{\mathbb P}^1$, the map
$\phi_\gamma\big\vert_{E_x}\,:\, E_x\, \longrightarrow\, E_{\gamma(x)}$
is a $\mathbb C$--linear isomorphism.

\begin{lemma}\label{lem1}
Let $E$ be a $G$--equivariant holomorphic vector bundle over ${\mathbb C}{\mathbb P}^1$. Let
$W\,\subset\, E$ be a $G$--invariant holomorphic subbundle. Assume that the short exact sequence
$$
0\, \longrightarrow\, W \, \longrightarrow\,
E \, \stackrel{q}{\longrightarrow}\, Q\,:=\, E/W \,\longrightarrow \, 0
$$
splits holomorphically, meaning 
there is a holomorphic subbundle $F \, \subset\, E$ such that the projection
$q\big\vert_F \, :\, F\,\longrightarrow\, Q$ is an isomorphism. Then $W$ has a $G$--invariant
direct summand (in other words, $F$ can be chosen to be preserved by the action of $G$ on $E$).
\end{lemma}

\begin{proof}
Since the action of $G$ on $E$ preserves $W$, the action of $G$ on $E$ produces an action of
$G$ on the quotient $Q$. The action of any $\gamma\, \in\, G$ on $Q$, given by $\phi_\gamma$ (see
\eqref{e3}) will be denoted by $\phi'_\gamma$. We have the holomorphic homomorphism
$$f \, \, :=\, \, (q\big\vert_F)^{-1} \,\, :\,\, Q \,\, \longrightarrow\,\, E ,$$
where $F$ is the subbundle in the lemma. Note that the equality $q\circ f\,=\, {\rm Id}_Q$ holds.

Consider the map $\phi$ in \eqref{e2} giving the $G$--equivariant structure of $E$. We have the homomorphism
$$
\widetilde{f}\,\,:\,\, Q \,\,\longrightarrow\,\, E,\ \ \, v\,\, \longmapsto\,\,
\frac{1}{\# G} \sum_{\gamma\in G} \phi^{-1}_\gamma\circ f\circ \phi'_\gamma.
$$
It is straightforward to check the following:
\begin{itemize}
\item $q\circ\widetilde{f}\,=\, {\rm Id}_Q$, and

\item the action of $G$ on $E$ preserves $\widetilde{f}(Q)$.
\end{itemize}
Therefore, the image $\widetilde{f}(Q)$ is a $G$--invariant direct summand of $W$.
\end{proof}

\begin{corollary}\label{cor1}
Let $E$ be a $G$--equivariant holomorphic vector bundle over ${\mathbb C}{\mathbb P}^1$ and
$$
0\, =\, W_0\, \subset\, W_1 \, \subset\, \cdots\, \subset\, E_{\ell-1}\, \subset\, E_\ell\,=\, E
$$
a filtration of $G$--equivariant subbundles of it, such that the
filtration splits holomorphically. Then the filtration splits holomorphically
$G$--equivariantly. In particular, $E$ is holomorphically $G$--equivariantly isomorphic to
the direct sum $\bigoplus_{i=1}^\ell W_i/W_{i-1}$.
\end{corollary}

\begin{proof}
Applying Lemma \ref{lem1} to the subbundle $W_{i-1}\, \subset\, W_i$, where $1\, \leq\, i\, \leq\, \ell$,
we conclude that $W_i\,=\, W_{i-1}\bigoplus (W_i/W_{i-1})$ as $G$--equivariant holomorphic vector bundles.
This immediately implies that the filtration in the statement of the corollary splits holomorphically
$G$--equivariantly.
\end{proof}

\section{Polystable equivariant bundles}

Any holomorphic vector bundle $V$ over ${\mathbb C}{\mathbb P}^1$ of rank $r$ decomposes into a direct sum of
the form $V\,=\,\bigoplus_{i=1}^r {\mathcal O}_{{\mathbb C}{\mathbb P}^1}(n_i)$, where $n_i$ are integers
\cite[p.~122, Th\'eor\`eme 1.1]{Gr}. We note that $V$ is semistable if $n_1\,=\, \cdots\,=\, n_r$. Also,
any semistable vector bundle on ${\mathbb C}{\mathbb P}^1$ is polystable. (See \cite[p.~14, Definition
1.2.12]{HL} and \cite[p.~23, Definition 1.5.4]{HL} for semistability and polystability respectively.)

\begin{proposition}\label{prop1}
Any $G$--equivariant holomorphic vector bundle $E$ on ${\mathbb C}{\mathbb P}^1$ is
a direct sum of $G$--equivariant polystable vector bundles.
\end{proposition}

\begin{proof}
If $E$ is semistable, then, as noted above, it is polystable. So we need to consider
the case where $E$ is not semistable. Let
\begin{equation}\label{ef}
0\, =\, W_0\, \subset\, W_1 \, \subset\, \cdots\, \subset\, E_{\ell-1}\, \subset\, E_\ell\,=\, E
\end{equation}
be the Harder--Narasimhan filtration of $E$ \cite[p.~16, Theorem 1.3.4]{HL} (see also \cite[p.~14,
Theorem 1.3.1]{HL}). Since $E$ is a direct sum of holomorphic line bundles, the filtration in \eqref{ef}
splits holomorphically. From the uniqueness property of the Harder--Narasimhan filtration it follows
immediately that each $W_i$ in \eqref{ef} is preserved by the action of $G$ on $E$. Therefore, from
Corollary \ref{cor1} we know that $E$ is holomorphically $G$--equivariantly isomorphic to
the direct sum $\bigoplus_{i=1}^\ell W_i/W_{i-1}$. Since $W_i/W_{i-1}$ is semistable, and any semistable
vector bundle over ${\mathbb C}{\mathbb P}^1$ is polystable, the proof is complete.
\end{proof}

\begin{theorem}\label{th1}
Set $G\,=\, {\mathbb Z}/n{\mathbb Z}$. Let $E$ be a $G$--equivariant holomorphic vector bundle over
${\mathbb C}{\mathbb P}^1$. Then $E$ splits into a direct sum of $G$--equivariant holomorphic line
bundles.

Any holomorphic vector bundle on ${\mathbb C}{\mathbb P}^1$ admits a $G$--equivariant structure.
\end{theorem}

\begin{proof}
Every holomorphic vector bundle over ${\mathbb C}{\mathbb P}^1$ decomposes into a direct sum of
holomorphic line bundles \cite[p.~122, Th\'eor\`eme 1.1]{Gr}. Moreover, any holomorphic line bundle
on ${\mathbb C}{\mathbb P}^1$ is
of the form ${\mathcal O}_{{\mathbb C}{\mathbb P}^1}(n)$ for some integer $n$. A $G$--equivariant
structure on a holomorphic vector bundle $E$ produces a $G$--equivariant structure on $E^*$ and
it also produces a $G$--equivariant structure on any tensor power of $E$. Therefore, to prove that
any holomorphic vector bundle $V$ on ${\mathbb C}{\mathbb P}^1$ admits a $G$--equivariant structure
it suffices to show that the two line bundles ${\mathcal O}_{{\mathbb C}{\mathbb P}^1}$ and
${\mathcal O}_{{\mathbb C}{\mathbb P}^1}(-1)$ admit a $G$--equivariant structure.

It was noted in Section \ref{se2} that $G$ is
conjugate to the subgroup of $\text{PGL}(2,{\mathbb C})$ generated by
$\begin{pmatrix}\alpha & 0\\ 0 & 1\end{pmatrix}$, where $\alpha$ is a primitive
$n$-th root of unity. So we assume that $G$ is the cyclic subgroup generated by
$\begin{pmatrix}\alpha & 0\\ 0 & 1\end{pmatrix}$, where $\alpha$ is a primitive
$n$-th root of unity.

The standard action of $\begin{pmatrix}\alpha & 0\\ 0 & 1\end{pmatrix}$
on ${\mathbb C}^2$ produces an action of $G$ on ${\mathcal O}_{{\mathbb C}{\mathbb P}^1}(-1)$.
Thus ${\mathcal O}_{{\mathbb C}{\mathbb P}^1}(-1)$ admit a $G$--equivariant structure. 

The trivial action of $\begin{pmatrix}\alpha & 0\\ 0 & 1\end{pmatrix}$ on $\mathbb C$ produces a
$G$--equivariant structure on ${\mathcal O}_{{\mathbb C}{\mathbb P}^1}$. Consequently,
every holomorphic vector bundle on ${\mathbb C}{\mathbb P}^1$ admits a $G$--equivariant structure.

Now we will show that any $G$--equivariant holomorphic vector bundle $E$ over
${\mathbb C}{\mathbb P}^1$ splits into a direct sum of $G$--equivariant holomorphic line
bundles. In view of Proposition \ref{prop1} we will assume that $E$ is polystable. Since the
holomorphic vector bundle $E$ decomposes into a direct sum of line bundles, we have
\begin{equation}\label{ed}
E\,\,=\,\, {\mathcal O}_{{\mathbb C}{\mathbb P}^1}(d)^{\oplus r}
\end{equation}
for some $d$, where $r\,=\, {\rm rank}(E)$. As observed above,
${\mathcal O}_{{\mathbb C}{\mathbb P}^1}(d)$ admits a $G$--equivariant structure. Fix
a $G$--equivariant structure on ${\mathcal O}_{{\mathbb C}{\mathbb P}^1}(d)$.

The $G$--equivariant structures on $E$ and ${\mathcal O}_{{\mathbb C}{\mathbb P}^1}(d)$ together
produce a $G$--equivariant structure on the vector bundle
$${\rm Hom}({\mathcal O}_{{\mathbb C}{\mathbb P}^1}(d), \, E)
\,=\, {\mathcal O}_{{\mathbb C}{\mathbb P}^1}(-d)\otimes E,$$
which, in turn, produces an action of $G$ on the vector space $H^0({\mathbb C}{\mathbb P}^1,\,
{\rm Hom}({\mathcal O}_{{\mathbb C}{\mathbb P}^1}(d), \, E))$. We have a natural homomorphism
\begin{equation}\label{eP}
\Phi\,\, :\,\, {\mathcal O}_{{\mathbb C}{\mathbb P}^1}(d)\otimes_{\mathbb C}
H^0({\mathbb C}{\mathbb P}^1,\,
{\rm Hom}({\mathcal O}_{{\mathbb C}{\mathbb P}^1}(d), \, E))\, \,\longrightarrow\,\, E
\end{equation}
that sends any $v\otimes s \, \in\, {\mathcal O}_{{\mathbb C}{\mathbb P}^1}(d)_x\otimes
H^0({\mathbb C}{\mathbb P}^1,\, {\rm Hom}({\mathcal O}_{{\mathbb C}{\mathbb P}^1}(d), \, E))$,
$x\, \in\,{\mathbb C}{\mathbb P}^1$, to $s(v)\, \in\, E_x$. From \eqref{ed} it follows immediately that
the holomorphic vector bundle ${\rm Hom}({\mathcal O}_{{\mathbb C}{\mathbb P}^1}(d), \, E)$ is
trivializable and $\Phi$ in \eqref{eP} is an isomorphism. Also, $\Phi$ is evidently $G$--equivariant.

Since $G$ is a finite cyclic group, the complex $G$--module 
$H^0({\mathbb C}{\mathbb P}^1,\, {\rm Hom}({\mathcal O}_{{\mathbb C}{\mathbb P}^1}(d), \, E))$
decomposes into a direct sum of one-dimensional complex $G$--modules. Let
$$
H^0({\mathbb C}{\mathbb P}^1,\, {\rm Hom}({\mathcal O}_{{\mathbb C}{\mathbb P}^1}(d), \, E))\,\,=\,\,
\bigoplus_{i=1}^N F_i,
$$
where each $F_i$ is a one-dimensional $G$--module and $N\,=\, \dim 
H^0({\mathbb C}{\mathbb P}^1,\, {\rm Hom}({\mathcal O}_{{\mathbb C}{\mathbb P}^1}(d), \, E))
\,=\, {\rm rank}(E)$. So
the isomorphism in \eqref{eP} gives an isomorphism
\begin{equation}\label{ed2}
\bigoplus_{i=1}^N {\mathcal O}_{{\mathbb C}{\mathbb P}^1}(d)\otimes_{\mathbb C} F_i\,\,
\stackrel{\sim}{\longrightarrow}\,\, E
\end{equation}
of $G$--equivariant holomorphic vector bundles. Note that \eqref{ed2} gives a holomorphic decomposition
of the $G$--equivariant vector bundle $E$ into a direct sum of $G$--equivariant holomorphic line bundles.
\end{proof}

\begin{remark}\label{rem1}
As in Theorem \ref{th1}, take $G\,=\, {\mathbb Z}/n{\mathbb Z}$. Let $L$ be a holomorphic line bundle
on ${\mathbb C}{\mathbb P}^1$. From Theorem \ref{th1} we know that $L$ admits a $G$--equivariant structure.
It is straightforward to check that any two $G$--equivariant structures on $L$ differ by multiplication
by a complex character of $G$. Indeed, this follows immediately from the fact that the group of holomorphic
automorphisms of $L$ is the multiplicative group ${\mathbb C}^\star\,=\, {\mathbb C}\setminus\{0\}$.
\end{remark}

\section{The Case of $G = ({\mathbb Z}/2{\mathbb Z})
\oplus ({\mathbb Z}/2{\mathbb Z})$}

In this section we assume that $G$ is the subgroup of $\text{PGL}(2,{\mathbb C})$ generated by
the two matrices $A_1$ and $A_2$ in \eqref{e1}. In particular, $G$ is isomorphic to
$({\mathbb Z}/2{\mathbb Z})\oplus ({\mathbb Z}/2{\mathbb Z})$. Let
\begin{equation}\label{f1}
p\,\,:\,\, \text{GL}(2,{\mathbb C})\,\, \longrightarrow\,\, 
\text{GL}(2,{\mathbb C})/{\mathbb C}^{\star} \,\,=\,\, \text{PGL}(2,{\mathbb C})
\end{equation}
be the natural quotient map, where ${\mathbb C}^\star \,=\, {\mathbb C}\setminus\{0\}$ is
the multiplicative group.

\begin{lemma}\label{lem2}\mbox{}
\begin{enumerate}
\item The holomorphic line bundle ${\mathcal O}_{{\mathbb C}{\mathbb P}^1}(2)$
admits a $G$--equivariant structure.

\item The holomorphic vector bundle ${\mathcal O}_{{\mathbb C}{\mathbb P}^1}(-1)
\oplus {\mathcal O}_{{\mathbb C}{\mathbb P}^1}(-1)$ admits a $G$--equivariant structure.

\item The holomorphic line bundle ${\mathcal O}_{{\mathbb C}{\mathbb P}^1}(-1)$ does
not admit any $G$--equivariant structure.
\end{enumerate}
\end{lemma}

\begin{proof}
The standard action of $G\, \subset\, \text{PGL}(2,{\mathbb C})$ on ${\mathbb C}{\mathbb P}^1$
has a natural lift to an action of $G$ on the holomorphic tangent bundle
$T{\mathbb C}{\mathbb P}^1$. Since $T{\mathbb C}{\mathbb P}^1\,=\, 
{\mathcal O}_{{\mathbb C}{\mathbb P}^1}(2)$, we conclude that
${\mathcal O}_{{\mathbb C}{\mathbb P}^1}(2)$ admits a $G$--equivariant structure.

Consider the rank two vector bundle $${\mathcal O}_{{\mathbb C}{\mathbb P}^1}(-1)
\oplus {\mathcal O}_{{\mathbb C}{\mathbb P}^1}(-1)\,\, =\,\, 
{\mathcal O}_{{\mathbb C}{\mathbb P}^1}(-1)\otimes_{\mathbb C} ({\mathbb C}^2)^*\,\,=\,\,
{\rm Hom}(\underline{\mathbb C}^2,\,{\mathcal O}_{{\mathbb C}{\mathbb P}^1}(-1)),
$$
where $\underline{\mathbb C}^2\,=\, {\mathbb C}{\mathbb P}^1\times {\mathbb C}^2
\, \longrightarrow\, {\mathbb C}{\mathbb P}^1$ is the trivial holomorphic vector bundle
of rank two on ${\mathbb C}{\mathbb P}^1$ with fiber ${\mathbb C}^2$.

Consider the standard actions of $\text{GL}(2,{\mathbb C})$ on
${\mathbb C}{\mathbb P}^1$ and ${\mathbb C}^2$. They together produce
a diagonal action of $\text{GL}(2,{\mathbb C})$ on
$\underline{\mathbb C}^2\,=\, {\mathbb C}{\mathbb P}^1\times {\mathbb C}^2$. 
 This action of $\text{GL}(2,{\mathbb C})$ on $\underline{\mathbb C}^2$ preserves
the tautological line subbundle
$$
{\mathcal O}_{{\mathbb C}{\mathbb P}^1}(-1)\,\, \subset\,\, \underline{\mathbb C}^2.
$$
the actions of $\text{GL}(2,{\mathbb C})$ on $\underline{\mathbb C}^2$ and
${\mathcal O}_{{\mathbb C}{\mathbb P}^1}(-1)$ together produce an action of
$\text{GL}(2,{\mathbb C})$ on 
\begin{equation}\label{f3}
{\rm Hom}(\underline{\mathbb C}^2,\,{\mathcal O}_{{\mathbb C}{\mathbb P}^1}(-1))\,\,=\,\,
{\mathcal O}_{{\mathbb C}{\mathbb P}^1}(-1)\oplus {\mathcal O}_{{\mathbb C}{\mathbb P}^1}(-1).
\end{equation}
Note that $\text{kernel}(p)\, \subset\, \text{GL}(2,{\mathbb C})$ acts trivially
on ${\rm Hom}(\underline{\mathbb C}^2,\,{\mathcal O}_{{\mathbb C}{\mathbb P}^1}(-1))$,
where $p$ is the projection in \eqref{f1}. Consequently, we get an action of
$\text{PGL}(2,{\mathbb C})\,=\, \text{GL}(2,{\mathbb C})/\text{kernel}(p)$ on
the vector bundle in \eqref{f3}. Therefore, the vector bundle in \eqref{f3} admits a
$G$--equivariant structure.

Now we will show that ${\mathcal O}_{{\mathbb C}{\mathbb P}^1}(-1)$ does
not admit a $G$--equivariant structure.

To prove by contradiction, assume that ${\mathcal O}_{{\mathbb C}{\mathbb P}^1}(-1)$
admits a $G$--equivariant structure. This implies that ${\mathcal O}_{{\mathbb C}{\mathbb P}^1}(1)$
admits a $G$--equivariant structure. Hence $G$ acts on
$$
H^0({\mathbb C}{\mathbb P}^1, \, {\mathcal O}_{{\mathbb C}{\mathbb P}^1}(1))
\,=\, {\mathbb C}^2.
$$
This action of $G$ on ${\mathbb C}^2$ makes $G$ a subgroup of $\text{GL}(2,{\mathbb C})$.

Any lift of $A_1\, :=\, \begin{pmatrix}-1 & 0\\ 0 & 1\end{pmatrix}$ (see \eqref{e1}) to
$\text{GL}(2,{\mathbb C})$ is of the form $A_{1,s}\,=\, \begin{pmatrix}-s & 0\\ 0 & s\end{pmatrix}$, with
$s\, \in\, {\mathbb C}^\star$. Also, any lift of $A_2\, :=\,\begin{pmatrix}0 & 1\\ 1 & 0\end{pmatrix}$
to $\text{GL}(2,{\mathbb C})$ is of the form 
$A_{2,t}\,:=\, \begin{pmatrix}0 & t\\ t & 0\end{pmatrix}$, with $t\, \in\, {\mathbb C}^\star$. Now 
$$
A_{1,s}A_{2,t} \,\,\not=\,\, A_{2,t}A_{1,s}
$$
for all $s,\, t$, and hence $G\,=\, ({\mathbb Z}/2{\mathbb Z})\oplus ({\mathbb Z}/2{\mathbb Z})$
can't lift to a subgroup of $\text{GL}(2,{\mathbb C})$. Consequently,
${\mathcal O}_{{\mathbb C}{\mathbb P}^1}(-1)$ does not admit any $G$--equivariant structure.
\end{proof}

\begin{theorem}\label{thm2}
Any $G$--equivariant vector bundle $E$ admits a holomorphic decomposition into $G$--equivariant
vector bundles
$$
E \,\,=\,\, \left(\bigoplus_{i=1}^M {\mathcal O}_{{\mathbb C}{\mathbb P}^1}(2m_i)\right)
\oplus \left(\bigoplus_{j=1}^N {\mathcal O}_{{\mathbb C}{\mathbb P}^1}(2n_i+1)^{\oplus 2}\right),
$$
where $m_i$ and $n_i$ are integers and $M+2N\,=\, {\rm rank}(E)$.
\end{theorem}

\begin{proof}
In view of Proposition \ref{prop1} we assume that $E$ is polystable. So the holomorphic vector
bundle $E$ has the following description:
\begin{equation}\label{f4}
E\,\,=\,\, {\mathcal O}_{{\mathbb C}{\mathbb P}^1}(b)^{\oplus r},
\end{equation}
where $b\, \in\, {\mathbb Z}$ and $r\,=\, {\rm rank}(E)$.

\textbf{Case 1:\, $b$ in \eqref{f4} is an even integer.}\, Assume that $b\,=\, 2m$. Since
${\mathcal O}_{{\mathbb C}{\mathbb P}^1}(2)$ admits a $G$--equivariant structure (see Lemma
\ref{lem2}(1)), we know that the holomorphic line bundle
\begin{equation}\label{f5}
L\,:=\, {\mathcal O}_{{\mathbb C}{\mathbb P}^1}(2m)\,=\, {\mathcal O}_{{\mathbb C}{\mathbb P}^1}(2)^{\otimes m}
\end{equation}
admits a $G$--equivariant structure; fix a $G$--equivariant structure on $L$ defined in \eqref{f5}.
The $G$--equivariant structures of $E$ and $L$ together define a $G$--equivariant structure on
$$
{\rm Hom}(L,\, E)\,\,=\,\, E\otimes L^*.
$$
Note that ${\rm Hom}(L,\, E)$ is a holomorphically trivializable vector bundle of rank $r$. Consider
$$
{\mathbb V}\,\, :=\,\, H^0({\mathbb C}{\mathbb P}^1,\, {\rm Hom}(L,\, E)).
$$
The action of $G$ on ${\rm Hom}(L,\, E)$ produces an action of $G$ on the vector space $\mathbb V$. As
in \eqref{eP}, let
\begin{equation}\label{f6}
\Phi\,\,:\,\, L\otimes_{\mathbb C} {\mathbb V}\,\, \longrightarrow\,\, E
\end{equation}
be the homomorphism that sends any $\ell\otimes s$, where $\ell\, \in\, L_x$ with $x\,\in\,
{\mathbb C}{\mathbb P}^1$ and $s\, \in\, {\mathbb V}$, to $s(\ell)\, \in\, E_x$.
We note that $\Phi$ in \eqref{f6} is an isomorphism. Also, $\Phi$ is $G$--equivariant.

Since $G$ is a finite abelian group, the $G$--module $\mathbb V$ decomposes into a direct sum of
$1$--dimensional complex $G$--modules. Fix a decomposition
$$
{\mathbb V}\,\,=\,\, \bigoplus_{i=1}^r F_i,
$$
where each $F_i$ is a one-dimensional complex $G$--module. Now
the isomorphism $\Phi$ in \eqref{f6} gives an isomorphism
$$
\Phi\,\, :\,\, \bigoplus_{i=1}^r L\otimes_{\mathbb C} F_i \,\, \longrightarrow\,\, E.
$$
So the $G$--equivariant vector bundle $E$ decomposes into a direct sum $G$--equivariant line bundles.

\textbf{Case 2:\, $b$ in \eqref{f4} is an odd integer.}\, Assume that $b\,=\, 2m+1$. Let
\begin{equation}\label{f7}
L\, \,:=\, \,{\mathcal O}_{{\mathbb C}{\mathbb P}^1}(b) \,\,=\,\, {\mathcal O}_{{\mathbb C}{\mathbb P}^1}(2m+1)
\end{equation}
be the holomorphic line bundle on ${\mathbb C}{\mathbb P}^1$.
Since ${\mathcal O}_{{\mathbb C}{\mathbb P}^1}(-2m-2)$ admits a $G$--equivariant structure
(see Lemma \ref{lem2}(1)) and ${\mathcal O}_{{\mathbb C}{\mathbb P}^1}(-1)\,=\,
{\mathcal O}_{{\mathbb C}{\mathbb P}^1}(-2m-2)\otimes L$ does not admit any $G$--equivariant structure
(see Lemma \ref{lem2}(3)), we know that $L$ in \eqref{f7} does not admit any $G$--equivariant structure.

Let 
$$
\widetilde{G}\,\, \subset\,\, \text{GL}(2, {\mathbb C})
$$
be the subgroup generated by the two elements $A_1,\, A_2\, \in\, \text{GL}(2, {\mathbb C})$ in
\eqref{e1}. So we have a central extension
\begin{equation}\label{f8}
0\, \longrightarrow\, {\mathbb Z}/2{\mathbb Z}\,=\, \pm {\rm I}\, \longrightarrow\,
\widetilde{G}\, \stackrel{q}{\longrightarrow}\, G \, \longrightarrow\, 0.
\end{equation}

The natural action of $\text{GL}(2,{\mathbb C})$ on ${\mathbb C}{\mathbb P}^1$ has a natural lift to
an action of $\text{GL}(2,{\mathbb C})$ on the line bundle ${\mathcal O}_{{\mathbb C}{\mathbb P}^1}(-1)$.
Therefore, $\widetilde{G}$ acts on the line bundle $L$ in \eqref{f7}. Using the projection $q$ in
\eqref{f8}, the action of $G$ on $E$ produces an action of $\widetilde G$ on $E$. The actions of
$\widetilde G$ on $L$ and $E$ together produce an action of $\widetilde G$ on 
$$
{\rm Hom}(L,\, E)\,\,=\,\, E\otimes L^*.
$$
This action of $\widetilde G$ on ${\rm Hom}(L,\, E)$ produces an action of $\widetilde G$ on
the vector space
\begin{equation}\label{f9}
{\mathbb V}\,\, :=\,\, H^0({\mathbb C}{\mathbb P}^1,\, {\rm Hom}(L,\, E)).
\end{equation}

Since $A^2_1\,=\, {\rm I}\,=\, A^2_2$ (see \eqref{e1}), the eigenvalues of $A_1$ and $A_2$ are
$\pm 1$. Let $v_1\, \in\, {\mathbb V}\setminus \{0\}$ be an eigenvector, for the
eigenvalue $\lambda$, for the action of $A_1\, \in\, \widetilde{G}$ on $\mathbb V$. Since
\begin{equation}\label{f10}
A_1A_2\,\,=\,\, - A_2A_1,
\end{equation}
we know $A_2(v_1)\, \notin\, {\mathbb C}\cdot v_1$. We have $A_1(A_2(v_1))\,=\, - A_2(A_1(v_1))$
(see \eqref{f10}), and hence $A_1(A_2(v_1))\,=\, - \lambda A_2(v_1)$. So $A_2(v_1)$ is an eigenvector,
for the eigenvalue $-\lambda$, for the action of $A_1\, \in\, \widetilde{G}$ on $\mathbb V$. Also,
note that $A^2_2(v_1)\,=\, v_1$. Consequently, the two-dimensional subspace
$$
{\mathbb V}_1\,\, :=\,\, {\mathbb C}\cdot v_1 \oplus {\mathbb C}\cdot A_2(v_1) \,\, \subset\,\, {\mathbb V}
$$
is preserved by the action of $\widetilde G$ on $\mathbb V$.

Next take another eigenvector
$v_2\, \in\, {\mathbb V}\setminus {\mathbb V}_1$ for the action of $A_1$ on $\mathbb V$ (note that
the action of $A_1$ is diagonalizable). Then we have
$$
{\mathbb V}_2\,\, :=\,\, {\mathbb C}\cdot v_2 \oplus {\mathbb C}\cdot A_2(v_2) \,\, \subset\,\, {\mathbb V}
$$
such that ${\widetilde G}\cdot {\mathbb V}_2\,=\, {\mathbb V}_2$
and ${\mathbb V}_1\oplus {\mathbb V}_2\, \subset\, {\mathbb V}$.

Proceeding this way, we get a decomposition of the $\widetilde G$--module $\mathbb V$
\begin{equation}\label{f11}
{\mathbb V}\,\,=\,\, \bigoplus_{j=1}^{r'} {\mathbb V}_j,
\end{equation}
such that $\dim {\mathbb V}_j\,=\, 2$ for all $1\,\leq\, j\, \leq\, r'$. Note that the
holomorphic vector bundle ${\rm Hom}(L,\, E)$ is trivializable, and hence $\dim {\mathbb V}\,=\,
{\rm rank}(E)\,=\, r$. Therefore, we have $2r'\,=\, r$ (see \eqref{f11}).

As before, we have an holomorphic isomorphism 
$$
\Phi\,\,:\,\, L\otimes_{\mathbb C} {\mathbb V}\,\, \longrightarrow\,\, E
$$
that sends any $\ell\otimes s$, where $\ell\, \in\, L_x$ with $x\,\in\,
{\mathbb C}{\mathbb P}^1$ and $s\, \in\, {\mathbb V}$, to $s(\ell)\, \in\, E_x$. This isomorphism
$\Phi$ is $\widetilde G$--equivariant. Using $\Phi$, the decomposition in \eqref{f11} produces
a decomposition
\begin{equation}\label{el}
E\,\,=\,\, \bigoplus_{j=1}^{r'} L\otimes_{\mathbb C} {\mathbb V}_j
\end{equation}
of the $\widetilde G$--equivariant vector bundle $E$. Since the action of $\widetilde G$ on $E$ factors
through the quotient $G$ in \eqref{f8}, the action of $\widetilde G$ on each $G$--equivariant subbundle
$L\otimes_{\mathbb C} {\mathbb V}_j\, \subset\, E$ in \eqref{el} factors through $G$.
This completes the proof.
\end{proof}

\begin{remark}\label{rem2}
For the same reason as in Remark \ref{rem1}, any two $G$--equivariant structures on 
${\mathcal O}_{{\mathbb C}{\mathbb P}^1}(2m)$ differ by multiplication by a complex character of $G$.
\end{remark}

\begin{remark}\label{rem3}
Note that the group of all holomorphic automorphisms of
${\mathcal O}_{{\mathbb C}{\mathbb P}^1}(2n+1)^{\oplus 2}$ is identified with $\text{GL}(2,{\mathbb C})$.
Let $\mathcal G$ denote the group of all pairs of the form $(g,\, \rho)$, where $g\, \in\, G$ and
$\rho\, :\, {\mathcal O}_{{\mathbb C}{\mathbb P}^1}(2n+1)^{\oplus 2}\, \longrightarrow\,
{\mathcal O}_{{\mathbb C}{\mathbb P}^1}(2n+1)^{\oplus 2}$ is a holomorphic isomorphism of vector bundles
over the automorphism of ${\mathbb C}{\mathbb P}^1$ given by $g$. So we get a short exact sequence of groups
$$
0\, \longrightarrow\, \text{GL}(2,{\mathbb C})\, \longrightarrow\, {\mathcal G}\,
\stackrel{\varphi}{\longrightarrow}\,G
\, \longrightarrow\, 0;
$$
the above projection $\varphi$ sends any $(g,\, \rho)\,\in\, {\mathcal G}$ to $g$. We note that
a $G$--equivariant structure on ${\mathcal O}_{{\mathbb C}{\mathbb P}^1}(2n+1)^{\oplus 2}$ is a
homomorphism $\eta\, :\, G\, \longrightarrow\, {\mathcal G}$ such that $\varphi\circ\eta\,=\, {\rm Id}_G$.
Let
$$
{\mathbb S}
$$
denote the space of all homomorphisms $\eta\, :\, G\, \longrightarrow\, {\mathcal G}$ such that
$\varphi\circ\eta\,=\, {\rm Id}_G$. The group $\text{GL}(2,{\mathbb C})$ acts on $\mathbb S$: The action
of any $A\, \in\, \text{GL}(2,{\mathbb C})$ sends any homomorphism $\eta\, \in\, {\mathbb S}$ to the
homomorphism $g\, \longmapsto\, A^{-1}\eta(g)A$. The corresponding quotient space
${\mathbb S}/\text{GL}(2,{\mathbb C})$ is identified with the space of all isomorphism classes of
$G$--equivariant structures on ${\mathcal O}_{{\mathbb C}{\mathbb P}^1}(2n+1)^{\oplus 2}$. 
\end{remark}

\section*{Acknowledgements}

The first author is partially supported by a J. C. Bose Fellowship (JBR/2023/000003).

\section*{Declarations}

On behalf of all authors, the corresponding author states that there is no conflict of interest.


\begin{thebibliography}{AAAA}

\bibitem[DF1]{DF1} C. De Concini and F. Fagnani, Symmetries of differential behaviors and
finite group actions on free modules over a polynomial ring, {\it Math. Control Signal Systems}
{\bf 6} (1993), 307--321.

\bibitem[DF2]{DF2} C. De Concini and F. Fagnani, Equivariant vector 
bundles over affine subsets of the projective line, \textit{Annali
Scuola Norm. Sup. Pisa} \textbf{22} (1995), 341--361.

\bibitem[Do]{Do} I. V. Dolgachev, {\it MacKay correspondence}, 
https://dept.math.lsa.umich.edu/$\sim$idolga/McKaybook.pdf.

\bibitem[FW]{FW} F. Fagnani and J. C. Willems, Representations of time-reversible systems,
{\it Jour. Mathematical Systems, Estimation and Control} {\bf 1} (1991), 5--28.

\bibitem[Gr]{Gr} A. Grothendieck, Sur la classification des
fibr\'es holomorphes sur la sph\`ere de Riemann, \textit{Amer.
Jour. Math.} \textbf{79} (1957), 121--138.

\bibitem[HL]{HL} D. Huybrechts and M. Lehn, {\it The geometry of moduli spaces of sheaves}, Aspects
of Mathematics, E31, Friedr. Vieweg~\&~Sohn, Braunschweig, 1997.

\bibitem[Mo]{Mo} L. Moser-Jauslin, Triviality of certain equivariant vector bundles for finite cyclic
groups, {\it C. R. Acad. Sci. Paris} {\bf 317} (1993), 139--144

\bibitem[Sp]{Sp} T. A. Springer, {\it Invariant theory}, Lecture Notes in Math., Vol. 585
Springer-Verlag, Berlin-New York, 1977.

\end{thebibliography}
\end{document}